\documentclass[notitlepage,leqno,11pt]{article}
\usepackage{latexsym}
\usepackage{amsmath}
\usepackage{amsfonts}
\usepackage{amssymb}
\usepackage{amscd}

\renewcommand{\d}{\delta }
\newcommand{\D }{\Delta }

\newcommand{\e }{\varepsilon }

\renewcommand{\i }{\iota}

\renewcommand{\l }{\lambda }

\newcommand{\n }{\nabla }

\renewcommand{\o }{\omega }

\newcommand{\ov}{\overline}
\newcommand{\intbar}{\mathop{\int\makebox(-13.5,0){\rule[4pt]{.7em}{0.3pt}}%
\kern-6pt}\nolimits}
\newcommand{\wtilde }{\widetilde}

\newcommand{\be}{\begin{equation}}
\newcommand{\ee}{\end{equation}}
\newcommand{\bes}{\begin{equation*}}
\newcommand{\ees}{\end{equation*}}
\newcommand{\ba}{\begin{eqnarray}}
\newcommand{\ea}{\end{eqnarray}}
\newcommand{\bas}{\begin{eqnarray*}}
\newcommand{\eas}{\end{eqnarray*}}
\newenvironment{pf}{\noindent{ \bf  Proof}.\enspace}{\rule{2mm}{2mm}\medskip}
\newenvironment{pfn}{\noindent{\bf  Proof}}{\rule{2mm}{2mm}\medskip}

\newcommand{\R}{\mathbb{R}}

\newcommand{\N}{\mathbb{N}}
\renewcommand{\o }{\omega }

\usepackage{color}

\begin{document}

\author{Mohameden AHMEDOU$^{a}$\thanks{  E-mail: \texttt{Mohameden.Ahmedou@math.uni-giessen.de}}, Mohamed BEN AYED$^{b,c}$\thanks{E-mails: \texttt{M.BenAyed@qu.edu.sa and Mohamed.Benayed@fss.rnu.tn}} and
 Khalil EL MEHDI$^{b,d}$\thanks{ E-mail : \texttt{K.Jiyid@qu.edu.sa} } \\
 {\footnotesize
a : Department of Mathematics, Giessen University, Arndtstrasse 2, 35392, Giessen, Germany.}\qquad\quad\quad\\
{\footnotesize
b : Department of Mathematics, College of Science, Qassim University, Buraydah 51452, Saudi Arabia.}\\
{\footnotesize
c :  Facult{\'e} des Sciences de Sfax, Universit\'e de Sfax, Route Soukra, Sfax, Tunisia.}\qquad\quad\qquad\qquad\quad\\
{\footnotesize
 d : Facult\'e des Sciences et Techniques, Universit\'e de Nouakchott, Nouakchott, Mauritania.}\quad \qquad \quad
\quad}

\date{}

\bigskip

\title{\bf  On the   Nirenberg problem on spheres:\\
Arbitrarily  many solutions in a perturbative setting}

\newtheorem{lem}{Lemma}[section]
\newtheorem{pro}[lem]{Proposition}
\newtheorem{thm}[lem]{Theorem}
\newtheorem{rem}[lem]{Remark}
\newtheorem{cor}[lem]{Corollary}
\newtheorem{df}[lem]{Definition}

\maketitle

\bigskip

\bigskip

\noindent{\bf Abstract:} Given a smooth positive function $K$ on the standard sphere $(\mathbb{S}^n,g_0)$, we use   Morse theoretical methods and counting index formulae  to prove that, under generic conditions  on the function $K$, there are  arbitrarily many metrics $g$ conformally  equivalent to   $g_0$ and whose scalar curvature is given by the function $K$  provided that the  function is sufficiently close to the scalar curvature of $g_0$. Our approach leverages a comprehensive characterization of blowing-up solutions of a subcritical approximation, along with various Morse relations involving their indices. Notably, this multiplicity result is achieved without relying on any symmetry or periodicity assumptions about the function $K$.\\
 \medskip


\tableofcontents

\section{Introduction and main Results}

In the academic year $1969-1970$ Louis Nirenberg, likely  motivated by the Yamabe probem posed  the following question: Given  a smooth positive function  $K$ on the standard  sphere $\mathbb{S}^n; n \geq 3$ endowed with its standard metric $g_0$, does there exist a riemannian metric $g$ conformally equivalent to $g_0$ such that the scalar curvature  $R_g$ with respect to $g$ is given by the function $K$? Such a geometric question has a nice PDE  interpretation. Indeed writing $g = u^{\frac{4}{n-2}} g_0$, the Nirenberg problem  amounts to  solving  the following nonlinear problem involving the critical Sobolev exponent:
$$
(\mathcal{N}_K) \qquad L_{g_0} u \, = \, K u^{(n+2)/(n-2)}, \quad u > 0 \quad\mbox{in}\quad \mathbb{S}^n,
$$
where  $ L_{g_{0}}:= -\D_{g_{0}}  \, + \, {n(n-2)}/{4}  $ denotes the conformal Laplacian.

Such a problem has been intensively studied in the last half century. See \cite{Au76, AH91, Bahri-Invariant, BC2, BCCH, BN83, CY, CGY, CL, CL1, CL2, Chen-Xu, Hebey, KW1, yyli1, yyli2, SZ} and the references therein.

Shortly after Nirenberg posed his question, Kazdan and Warner identified topological obstructions  \cite{BEZ1, KW1}, showing that the Nirenberg problem cannot be solved for every function $K$. Therefore, the central inquiry revolves around identifying sufficient conditions on $K$ under which the problem is solvable. To provide conditions on the function  $K$, under which the Nirenberg problem is solvable, some   Euler-Poincar\'e type criteria have been established. Notable results include those  of Bahri-Coron \cite{BC2}, Chang-Gursky-Yang \cite{CGY} on $\mathbb{S}^3$, Ben Ayed et al \cite{BCCH} and Yanyan Li \cite{yyli2} on $\mathbb{S}^4$.
Further existence results have been obtained for higher-dimensional spheres under various conditions: flatness near the critical points of $K$ \cite{yyli1}, in the perturbative setting \cite{CY}, and under pinching conditions \cite{Chen-Xu, Malchiodi-Mayer}.

To understand the analytical challenges posed by the Nirenberg problem, one must interpret its solutions as critical points of an Euler-Lagrange functional. Specifically, positive critical points of the functional
 \be\label{functional}
I_{K}(u) := \frac{1}{2}  \| u \|^2 - \frac{ n-2 }{ 2n}  \int_{ \mathbb{S}^n} K  | u |^{2n / (n-2 ) } ,
\ee 
defined on  $H^1(\mathbb{S}^n)$  equipped with the norm:
$$
\| u \|^2= \int_{\mathbb{S}^n} \mid \n u\mid^2 +\frac{n(n-2)}{4}
\int_{\mathbb{S}^n}  u^2,
$$ 
correspond to solutions of the Nirenberg problem. However, due to the presence of the Sobolev critical exponent, this functional does not satisfy the Palais-Smale condition. This means that there exists a sequence $(u_k)$ along which $I_K(u_k)$  is bounded, $\n I_K(u_k)$  goes to zero and $(u_k)$ does not converge. The reason for such a lack of compactness is the existence of almost solutions of the equation $(\mathcal{N}_K)$ \cite{Struwe}. 
 
One  way to overcome the lack of compactness of the functional $I_K$, which traces  back to Yamabe \cite{Yamabe}, is to lower the exponent and first consider the subcritical approximation
 $$
(\mathcal{N}_{K,\tau}) \qquad L_{g_0} u \, = \, K u^{((n+2)/(n-2))-\tau}, \quad u > 0 \quad\mbox{on}\quad \mathbb{S}^n , 
$$
and its associated Euler-Lagrange functional
$$
 I_{K,\tau}(u) := \frac{1}{2}  \| u \|^2 - \frac{ n-2 }{ 2n - \tau (n-2) }  \int_{ \mathbb{S}^n} K  | u |^{\frac{2n}{n-2}-\tau} , \quad u \in H^1(\mathbb{S}^n) . 
 $$
 
Using  elliptic estimates, one obtains that  either the solution $u_{\tau}$ remains uniformly bounded as   $\tau \to 0$ and hence converges strongly to a solution $u_{\infty}$ of the Nirenberg problem $(\mathcal{N}_K)$ or it blows up. In the latter  case,  following  the works of Schoen \cite{Schoen}, Yanyan Li \cite{yyli1, yyli2}, Chen-lin \cite{CL1,CL2} or Druet-Hebey-Robert \cite{DHR}, one performs a refined blow up analysis. For  energy-bounded solutions of $(\mathcal{N}_{K,\tau})$, two  scenarios can arise: Either the solution  $u_{\tau} $ converges weakly to zero or it converges weakly to $\o \neq 0$ which is a solution of $(\mathcal{N}_K)$.
 The case of a zero weak limit has been extensively studied, and the complete blow-up picture on high-dimensional spheres has been fully described only in recent papers by Malchiodi and Mayer \cite{MM, MM19}. Recently, we investigated in \cite{Blow-up sphere} the scenario where the sequence of blowing-up solutions exhibits a non-zero weak limit.
 
We point out that the blow up picture described in \cite{Blow-up sphere} contrasts with the Yamabe case where the function $K$ is a constant, say $K \equiv 1$. Indeed the level sets of the Yamabe functional  $I_1$ on the sphere are all contractible, the  critical points of $I_1$ are minima and $I_1$ does not have any \emph{critical point at infinity}, which are  non compact orbits of the gradient flow. This  observation is   quite striking if one considers the case where the function $K$ is close to one. Indeed thinking of a homotopy between $K$ and constant $1$ one observes that all the blowing up solutions constructed in Theorem  $1.2$ of \cite{Blow-up sphere} will disappear when $K \equiv 1$ without giving rise to \emph{critical points at infinity} for $I_1$. Such a situation may only occur if there  is a \emph{cancellation phenomenon} with genuine solutions of $(\mathcal{N}_K)$ and our goal in this paper is to find a class of functions  $K$ for which such solutions exist.
 Namely we consider the following class of functions: For $\e >0$, we define
 $$
 \mathcal{K}_\e= \Big\{K:\mathbb{S}^n\to \mathbb{R}\, : \, \, 0<K\in C^3(\mathbb{S}^n)\quad\mbox{and}\quad   \Big   ( \frac{K_{max}}{K_{min} }-1   \Big) <\e\Big\},
 $$
 where $ K_{max} := \max_{\mathbb{S}^n} K$ and $K_{min} := \min_{\mathbb{S}^n} K $.\\
 Furthermore we assume   $K \in \mathcal{K}_{\e}$ to  satisfy the following generic conditions:\\
$\textbf{(H1)}$ The critical  points $y_1$,..., $y_{m_1}$  of $K$ are non degenerate and  satisfy $\D K(y_i) \neq 0$.\\
$\textbf{(H2)}$ The critical points of the associated variational function $I_K$  are non degenerate. \\
$\textbf{(H3)}$ $m:=\# \mathcal{K}^{\infty} \geq 2$, where $\mathcal{K}^{\infty} $  is defined  by
\be\label{Kinfty} \mathcal{K}^\infty  := \{ y \in \mathbb{S}^n : \n K(y) =0 \mbox{ and } \D K(y) < 0\} \ee
 and $\#\mathcal{K}^{\infty} $ denotes the cardinal of the set $\mathcal{K}^{\infty} $.
 
First, we consider the case where the Hopf-Poincar\'e index formula related to the topological contribution of critical points of $I_K$ is equal to $1$, that is, 
$$ \text{Index}_K := \sum_{y \in \mathcal{K}^\infty} (-1)^{n - morse\, index (K,y)} = 1.
$$
We remark that, in this case, we have necessarily that $m:=\# \mathcal{K}^{\infty} \geq 3$ and $m$ has to be odd.

 Given an arbitrary integer $N$, we claim that, under appropriate conditions on $K$, for $\e$ small and $1\leq k\leq [\frac{N}{2}]$, problem
$(\mathcal{N}_K)$ has solutions whose energies are close to $\frac{2k}{n}S_n$, where 
\be \label{eq:sn} S_n= \left(n(n-2)\right)^{ n / 2} \int_{\R^n} \frac{1}{(1+ | x | ^2 )^n } dx \ee 
 is the Sobolev universal constant. Here 
 $(\mathcal{N}_K)$ is said to have  a solution $ \o $ whose energy is   close to $ \frac{2k}{n}S_n$ if we have 
$$ I_K ( \o ) \in [( 2k\frac{S_n}{n}(1-\eta)), ( 2k\frac{S_n}{n}(1+\eta)) ] $$
for some small positive constant $\eta$.\\
More precisely, our first multiplicity result reads as follows: 
\begin{thm}\label{th:t1}
Let $n \geq 7$ and $N \in \mathbb{N}$ be an integer. There exists $\varepsilon_N > 0$ such that, if $ 0 < \varepsilon < \varepsilon_N$  and $K \in \mathcal{K}_\varepsilon$ satisfying  $\textbf{(H1-H2-H3)}$ and 
$$  \text{Index}_K := \sum_{y \in \mathcal{K}^\infty} (-1)^{n - morse\, index (K,y)} = 1,$$
 then the following facts hold: \\
    for any $k\in\{1,2, \cdots ,[\frac{N}{2}]\}$, $(\mathcal{N}_K)$ has at least $ C_{k + \ell -1 } ^k $ solutions $ u_{j,k}$'s whose energies are close to $ \frac{2k}{n}S_n$. Here, $ \ell $ denotes $ \ell:= (m-1) / 2 $. \\
    In particular,  $(\mathcal{N}_K)$ has at least $  C_{ \ell + [N/2 ] } ^ \ell - 1  $ solutions.
\end{thm}

Next, we deal with the case where the Hopf-Poincar\'e index formula related to the topological contribution of critical points of $I_K$ is different from 1, that is, $  \text{Index}_K  \neq 1$.  We take an arbitrary integer $N$ and our goal is to show that, for $\e$ small and $1\leq k\leq N$, problem
$(\mathcal{N}_K)$ has solutions whose energies are close to $\frac{k}{n}S_n$. 
 Namely, we prove:

\begin{thm}\label{th:tt1}
Let $n \geq 7$ and $N \in \mathbb{N}$ be an integer. There exists $\varepsilon_N > 0$ such that, if $ 0 < \varepsilon < \varepsilon_N$  and $K \in \mathcal{K}_\varepsilon$ satisfying  $\textbf{(H1-H2-H3)}$ and 
$$  \text{Index}_K := \sum_{y \in \mathcal{K}^\infty} (-1)^{n - morse\, index (K,y)} \neq 1,$$
 then the following facts hold: \\
    for any $k\in\{1,\cdots, N \}$, $(\mathcal{N}_K)$ has at least a solution $ u_{k}$ whose energy is close to $ \frac{k}{n}S_n$. \\
    In particular,  $(\mathcal{N}_K)$ has at least $  N  $ solutions.
\end{thm}

\begin{rem} \begin{itemize}
\item  For each $k\in\{1,\cdots, N \}$, we can be more precise for the number of solutions whose energies are close to $ \frac{k}{n}S_n$. See Theorems \ref{th:t2} and \ref{th:t2bis} for precise statements.
\item If  $\# \mathcal{K}^{\infty}$ is even then one can deduce that   $ \text{Index}_K := \sum_{y \in \mathcal{K}^\infty} (-1)^{\i (y)} \neq 1$.
\end{itemize}
\end{rem}

Up to the knowledge of the authors the  results  presented above are the first multiplicity results for the Nirenberg problem  which do not rely on  periodicity or symmetry assumptions for the function $K$. 

The underlying idea behind our multiplicity results  can be  described as follows: If the total index of the blowing up solutions is different from one then one  derives, under assumption that  the function $K$ is close to one,   the existence of a solution whose energy level is close to the first Yamabe level $S_n/n$. Such a solution gives rise by Theorem $2$ of \cite{Blow-up sphere} to a blowing up solution near  the second Yamabe energy level $2 S_n/n$. Taking into account the topological contribution of all blowing  up solutions and the fact that for $K$ close enough to $1$ the nearby levels are topologically trivial one derives the existence of a second solution which again gives rise to a  blowing up solution at the third Yamabe level $3S_n/n$ and so one.  Using Morse equalities  relating to  the contribution of the solutions to  the difference of topology between the level sets of the associated Euler-Lagrange functional, one proves  the existence of  a new solution  of the Nirenberg problem  near every Yamabe level and hence a new blowing up solution of the approximate problem. In case where the above mentioned sum is equal to one a similar argument  provides the existence of a solution   near the second Yamabe level.

\bigskip

 The organization of the sequel is as follows :   in  Section $2$ we explore  the topology of certain  sublevels of the Euler-Lagrange functional $I_K$. In  Section $3$ we present  the proof of Theorem \ref{th:tt1} in the  special case while the contribution of the critical points of $K$ in the difference of topology between the level sets of  $I_K$ is spherical. Section $4$ is devoted to derive several  Morse inequalities and index counting  formulae relating the critical points of $I_K$ via  their Morse indices. Finally in    Section $5$ we provide the full  proofs of our main  multiplicity results.

\section{Topology of some level sets of the Euler-Lagrange functional}

 We first,  observe that, if $ u $ is a solution of $ ( \mathcal{N}_K ) $ then $ \beta ^{ (2-n) /4 } u $ is a solution of $ (  \mathcal{N}_{ \beta K } ) $, where $\beta >0$. Therefore, without loss of generality,  we can assume that $ K_{min } = 1 $  in the sequel. Therefore if  $ K $ in $ \mathcal{ K _\varepsilon }$, it follows that $ K $ is close to $ 1 $. Moreover since the  free functional $I_{K,\tau}$ is not lower bounded,
it is more convenient to consider the following functional
$$ J_{K,\tau} (u) := \frac{1}{ \Big( \int _{ \mathbb{S}^n} K | u | ^{ \frac{2n}{n-2} - \tau} \Big)^{2 (n-2) / ( 2n - \tau (n-2) )}} \qquad u \in \Sigma := \{ u\in H^1( \mathbb{S}^n) : \| u \| = 1 \}.$$
By easy computations, we can check that
\be \label{JKIK1} \n J_{K,\tau} (u) = 2 J_{K,\tau} (u)^{ 1 - \frac{ 2n - \tau (n-2) }{ 8 - 2 \tau (n-2)}} \n I_{K,\tau} \Big( J_{K,\tau}(u) ^{\frac{ 2 n - \tau (n-2) }{ 8 - 2 \tau (n-2 )}} u \Big) \qquad \forall \, \, u \in \Sigma. \ee
Thus, up to a multiplicative constant, the positive critical points of $J_{K,\tau} $  are in one to one correspondence with  the positive critical points of $I_{K,\tau} $.
Furthermore, for $ \tau = 0 $, if $ u $ is a positive critical point of $ J_K $, it follows that $ \o := J_K(u) ^{n/4} u $ is a  positive critical point of $ I_K $ (that is a solution of $ (\mathcal{N} _K ) $) and we have 
\be \label{JKIK2} I_K ( \o ) = \frac{1}{n} J_K(u)^{n/2} \qquad \mbox{ with }  \o := J_K(u) ^{n/4} u . \ee 

Our aim in this section is to obtain some useful information about some sublevel sets of the functional $J_{K,\tau}$ under the assumption that the function $K$ is close to a constant. 
First,   we recall the following   abstract lemma due to Malchiodi-Mayer  \cite{Malchiodi-Mayer} (Proposition $3.1$), which reads as follows:

\begin{lem}[\cite{Malchiodi-Mayer}]\label{deform}
Let $\underline{A}$ and $\ov{A}:=  ({ K_{\max}}/{ K_{\min}})^{(n-2)/n} \, \underline{A}$. Assume that $J_{K,\tau}$ does not have any critical point  in the set $J_{K,\tau}^{\ov{A}} \setminus J_{K,\tau}^{\underline{A}}$ where $J_{K,\tau}^A := \{ u: J_{K,\tau}(u) < A\}$. Then for each $c\in [ \underline{A},\ov{A}]$, the level set  $J_{K,\tau}^c$ is contractible.
\end{lem}
The crucial ingredient in the proof of our multiplicity Theorems is the following result:
\begin{pro}\label{l:C} Let $N\in \mathbb{N}$  and $\eta$ be a positive constant satisfying $ \eta < 1/(2N+1)$. There exists $\e_{N,\eta} > 0$ such that, for
 $$  0 <\e <\min \Big(\e_{N,\eta}, \left(\frac{N+1}{N}\right)^{2/(n-2)}\left(\frac{1-\eta}{1+\eta}\right)^{2/(n-2)} - 1 \Big)$$  and  $K\in \mathcal{K}_\e$,  the sublevel sets $J_{K, \tau}^C$, for $\tau $ small enough,  are  contractible for any
$$
C\in \left[ ( k S_n (1+\eta) ) ^{2/n} , ( (k+1) S_n (1-\eta) ) ^{2/n}\right] \quad \forall \, \, k\in\{1,...,N\}, $$ where $S_n$ is  the Sobolev constant defined in  \eqref{eq:sn}.
\end{pro}
\begin{pf}
It follows from the quantization  of Palais-Smale sequences  for the Yamabe functional  $J_1$ (see \cite{Lions, Struwe}) and the classification of  Yamabe solutions on the sphere (see  \cite{Obata, CGS})  that each Palais-Smale sequence $(u_q)_q$ (under the level $(N+1)^{2/n}S_n ^{2/n}$)
has to be close to a sum of $\tilde\d_{a_i, \lambda_i} $ and therefore $J_1(u_q)$ will be close to $k^{2/n}\times S_n^{2/n}$  for some $k \leq N$.
Here, for $ a \in \mathbb{S}^n$ and $ \lambda > 0 $, the function $\tilde\d_{a,\lambda}$ denotes the solution of the Yamabe Problem on $ \mathbb{S}^n$ and it is defined by 
\begin{equation}\label{eq:deltatilde}
\wtilde{\d}_{(a,\l)}(x)=c_0
\frac{\l^{(n-2)/{2}}}{\bigl(2+(\l^2-1)(1-\cos d_{g_0}(x,a))\bigr)^{(n-2)/{2}}}, \quad \mbox{ with }\, \,  c_0 := [ n(n-2)^{(n-2)/4} 
\end{equation} 
where $ d_0 $ is the geodesic distance on $\mathbb{S}^n$.

Hence $ | \n J_1(u) | \geq 3 C_{N,\eta}$ for some positive constant $ C_{N,\eta}$ and for each $u$ such that  $$
J_1(u) \in \bigcup_{1\leq k \leq N} \left[ ( k S_n (1+\eta) ) ^{2/n} , ( (k+1) S_n (1-\eta) ) ^{2/n}\right].
 $$

 Next, since $ K = 1 + O(\e) $, we derive that 
\begin{align*}
& J_K (u ) = ( 1+ O(\e)) J_1 (u) , \\
 & \langle \n  J_K (u )  , h \rangle = 2 J_K (u ) \left( \langle u , h \rangle - J_K (u ) ^{n/(n-2) } \int_{\mathbb{S}^n} K | u | ^{ (n+2)/(n-2) } h \right)
\end{align*}
which implies that $\n J_{K} \to \n J_1$ as $\e\to 0$. Thus,  there exists $\e_{N,\eta}$ such that for $0<\e\leq \e_{N,\eta}$ and $K\in \mathcal{K}_\e$, we have
\be \label{smale1} \vert \n J_{K}(u)\vert \geq 2 C_{N,\eta} \quad 
\mbox{ for }  u \mbox{  satisfying }
 J_{K}(u)\in \cup_{k=1 }^N \left[ ( k S_n (1+\eta) ) ^\frac{2}{n} , ( (k+1) S_n (1-\eta) ) ^\frac{2}{n}\right]   .\ee
Now, we need to prove that an analogue result is true for $ J_{K, \tau}$. To do so, let $ (u_q)_q$ be a positive Palais-Smale sequence for $ J_{K, \tau}$ under the level  $(N+1)^{2/n}S_n ^{2/n}$. On one hand, since $ J_{K, \tau}$ satisfies the Palais-Smale condition, it follows that $u_q$ has to be close to a positive critical point $ \o^  \tau $ of $ J_{K, \tau}$ and therefore we have 
\be\label{smale2}
 J_{K, \tau} (u_q) =  J_{K, \tau} (w ^ \tau) + o_q(1) \ee 
 where $ o_q (1) $ is a quantity which goes to zero as $ q \to \infty$.   Hence, $ ( \o^ \tau )_\tau $ is  a family of critical points of $ J_{K, \tau} $ under the level  $(N+1)^{2/n}S_n ^{2/n}$.  
 
 On the other hand, combining Theorem $1$ in \cite{MM19} and Theorem $1.1$ in \cite{Blow-up sphere}, we have that $\o ^ \tau$ either  converges strongly  to a positive critical point $\o_p$ of $ J_K $ and therefore we deduce that 
 $ J_{K, \tau} (w ^ \tau) = J_{K } (w _p) + o_\tau (1) $,   or it converges weakly to a positive critical point $\o_q$ of $ J_K $   (with $q < p$) and blows up at $p-q$ points $y_i$'s in $ \mathcal{K}^\infty$ and in this case we deduce that $ J_{K, \tau} (w ^ \tau) = J_{K } (w _q) + (p-q) S_n + o_\tau (1) $. \\
Here and in the sequel $\o_k$ denoted a solution of  $( \mathcal{N}_K )$ having an energy $J_{K}(\o_k)$ lying  in the interval  $[ ( k S_n(1 - \eta))^{2/n}, ( k S_n (1 + \eta) ) ^{2/n} ]$ for $k \in \{1, \cdots,N \}$.\\
Thus, we deduce that the energy of each Palais-Smale sequence of $ J_{K,\tau}$ has to be close to $ ( p S_n)^{2/n}$ for some $ p \in \{ 1, \cdots, N \}$. Therefore, for small $\e$ and small $ \tau $, 
$$ \vert \n J_{K, \tau}(u)\vert \geq  C_{N,\eta} \quad 
\mbox{ for }  u \mbox{  satisfying }
 J_{K, \tau}(u)\in \cup_{k=1 }^N \left[ ( k S_n (1+\eta) ) ^\frac{2}{n} , ( (k+1) S_n (1-\eta) ) ^\frac{2}{n}\right]   .$$

But, since $ J_{K,\tau}$ satisfies the Palais-Smale condition, it follows that $ J_{K,\tau}^{ ( (k+1) S_n (1-\eta) ) ^{2/n} }$ retracts by deformation onto $ J_{K,\tau}^{ ( k S_n (1+\eta) ) ^{2/n} }$. Lastly, to apply Lemma   \ref{deform}, we take $\underline{A} = ( k S_n (1+\eta) ) ^{2/n} $ and it is sufficient to get  $ \ov{A} \leq ( (k+1) S_n (1-\eta) ) ^{2/n} $ which implies that 
$$
\Big(\frac{k+1}{k} \Big)^{2/(n-2)} \Big(\frac{1-\eta}{ 1 + \eta} \Big)^{2/(n-2)} - 1 \geq \e  \quad \forall \, \, \, k \in \{ 1, \cdots, N\} .
$$
Since $ \frac{k+1}{k} \geq \frac{N+1}{N}$ for $ k \leq N$, our proposition follows.
\end{pf}

\section{  Proof of Theorem \ref{th:tt1}  in a simplified framework}

 The proofs of our results are basically based on Morse theoretical arguments and some counting Morse index of blowing up solutions of  the subcritical approximation. We notice that, when $ m:= \# \mathcal{K}^{\infty}= 2$, we clearly have $ \text{Index}_K \neq 1$ and therefore we are in the situation of Theorems \ref{th:t2} and \ref{th:t2bis} stated in Section $5$. Now, we are going to fix 
 some notation which will be used  in the proof of the theorems.  For $y\in \mathcal{K}^\infty$, we denote by
$$  \i ( y) = n - \mbox{ morse index } (K,y). $$
 Suppose that $\mathcal{K}^\infty = \{y_1, \cdots ,y_m\}$ with $y_1$ is a  global  maximum point of $K$ and  for $N \in \N$ and $\eta < \frac{1}{2N +1}$ we choose $\e$ as  in Proposition \ref{l:C}.
We denote by $u_p^\tau$ a solution of $( \mathcal{N}_{K,\tau})$ having an energy  $J_{K,\tau}(u_p ^ {\tau})$    in the interval   $[(p S_n(1 - \eta))^{2/n}, ( p S_n(1 + \eta))^{2/n}] $. Combining Theorem $1$ in \cite{MM19} and Theorem $1.1$ in \cite{Blow-up sphere}, we have that $u_p^\tau$ either  converges strongly  to a solution $\o_p$ of $( \mathcal{N}_K )$ or it converges weakly to a solution $\o_q$ of $(\mathcal{N}_K )$   (with $q < p$) and blows up at $p-q$ points $y_i$'s in $ \mathcal{K}^\infty$. \\
Here and in the sequel $\o_k$ denoted a solution of  $( \mathcal{N}_K )$ having an energy $J_{K}(\o_k)$ lying  in the interval  $[ ( k S_n(1 - \eta))^{2/n}, ( k S_n (1 + \eta) ) ^{2/n} ]$ for $k \in \{1, \cdots,N \}$.\\
Next observe that  if $u_p^\tau$ converges  strongly as $\tau \to 0$   to $\o_p$, then 
\be\label{KB1} m( u_p^\tau ) := \mbox{morse index } (J_{K,\tau}, u_p^\tau ) = \mbox{morse index } (J_K , \o_p ) := m (\o_p),
\ee
and when  $u_p^\tau$  converges weakly to $\o_q$ and blows up at $y_{i_{1}}, \cdots, y_{i_{(p-q)}}$ (in this case, we will denote it  by  $u_{\infty, p}^\tau$), we have that  $$ m( u_{\infty, p}^\tau ) = (p-q) + m(\o_q) + \sum_{j=1}^{p-q} \i (y_{i_j}).$$
However, if $u_p^\tau$  converges weakly to $ 0 $ and blows up at $y_{i_{1}}, \cdots, y_{i_{p}}$, we have that  
$$ m( u_{\infty, p}^\tau ) = (p - 1) + \sum_{j=1}^{p} \i (y_{i_j}).$$

Next we denote by 
\be\label{****}
\mathcal{S}^p=\{u_p^\tau\,:J_{K, \tau}(u_p^{\tau})\in \big [ ( p  S_n (1 - \eta) ) ^{2/n} , ( p S_n (1 + \eta) ) ^{2/n} \big]\,\,\mbox{and}\,\,u_p^\tau\to  \o_p\},
\ee
where $\o_p$ is a critical point of $ J _K := J_{K, 0} $.

 We set
\begin{equation}\label{mup}
   \mu_p :=
   \begin{cases}
  0, & \mbox{if }  \mathcal{S}^p = \emptyset ,  \\
   \sum_{ u_{ p }^\tau \in \mathcal{S}^ p } (-1)^{ m( u_{p}^\tau )}, & \mbox{otherwise}.
\end{cases}
\end{equation}
We observe that, when $\mathcal{S}^p\neq \emptyset$,  \eqref{KB1} implies that 
\begin{equation}\label{mupw}
\mu_p = \sum_{\o_p} (-1)^{m(\o_p)}.
\end{equation}
Hence, for some positive small constant $ \eta $, letting 
\be \label{critJK}
\text{Crit}^p(J_K) := \{ \o _p : \n J_K( \o_p) = 0 \text{ and } J_K(\o_p) \in  \big [ ( p  S_n (1 - \eta) ) ^{2/n} , ( p S_n (1 + \eta) ) ^{2/n} \big] \},
\ee
 we deduce  from \eqref{mupw} that 
\be\label{cardinalp}
\# ( \text{Crit}^p(J_K) ) \geq | \mu_p  | ,
\ee
where $ \#(F)$ denotes the cardinal of the set $ F $.\\
Therefore, to prove the theorems, we need to get some information about the values of the $ \mu_p $.
 
 \medskip
 
\begin{pfn}{ \bf   of Theorem \ref{th:tt1}  when $ m:= \# \mathcal{K}^{\infty}= 2$} \\
 For the sake of clarity we first prove the theorems in a simplified framework where the topology induced  by the function $K$  to the difference  of topology between the sublevels $J_{K}^{  ( S_n (1+ \eta) ) ^{2/n} }$  and $J_{K}^{ ( S_n (1 -\eta) ) ^{2/n} }$
   is spherical. Namely we consider the case where
 $\mathcal{K}^\infty = \{y_1,y_2\}$ with $y_1$ is a  global  maximum point of $K$  and $y_2$ is a critical point of $K$ of co-index $ \i_2:= n - morse(K,y_2)$.  \\
 To prove that $\mu_1 \neq 0$, we observe that, according to Proposition \ref{l:C}, we have that the level set  $J_{K,\tau}^{  (S_n (1+ \eta))^{2/n} }$ is contractible. Furthermore according to  Theorem 1 in \cite{MM} there are exactly two blowing up solutions  converging weakly to $0$ and  whose energy levels lie under this  level set. These solutions blow up   respectively at $y_1$ and $y_2$ and their Morse indices are respectively $0$ and $\i_2:= \i(y_2) $. Hence it follows from Euler-Poincar\'e Theorem that
 $$
 1 \, = \chi( J_{K,\tau}^{  (S_n (1+ \eta) ) ^{2/n} })\, =  \, \mu_1 \, +  1 + (-1)^{\i_2} .
 $$

 Therefore
 \be\label{muu1}
 \mu_1 \, = \,  -(-1)^{\i_2} \neq 0 .
 \ee
Hence the problem $(\mathcal{N}_K)$ has at least one solution under the  energy level $ S_n ^{2/n}(1+ \eta)^{2/n} $. We denote by $(\o_1^j)_{j=1}^{p_1}$ these solutions. \\
Again using Proposition \ref{l:C} we have that the level sets $J_{K,\tau}^{ ( 2 S_n (1+ \eta) ) ^{2/n} }$  and $J_{K,\tau}^{ ( 2 S_n (1 -\eta) ) ^{2/n} }$ are contractible. Furthermore    according to Theorem 1 in \cite{MM} and Theorem $1.2$ in \cite{Blow-up sphere}, blowing up solutions whose energy levels lie in the interval $[( 2 S_n (1- \eta))^{2/n}, (2 S_n (1+ \eta)) ^{2/n} ]$ are :
\begin{itemize}
\item the solution which converges weakly to $0$ and blows up at $y_1$ and $y_2$,
\item the solutions which converge weakly to a solution $\o_1^j$, $j=1,\cdots,p_1$, and blow up at $y_1$,
\item the solutions which converge weakly to a solution $\o_1^j$, $j=1,\cdots,p_1$, and blow up at $y_2$.
\end{itemize}
In addition, there are also the solutions which converge strongly to a solution $\o_2^j$.
 Therefore it follows from the Euler-Poincar\'e theorem that
$$
  1= \, \mu_2 \, + \, (-1)^{1+\i_2} \, + \, \sum_{j=1}^{p_1} (-1)^{1+ morse(\o_1^j)}  \,+  \,  \sum_{j=1}^{p_1} (-1)^{1+ \i_2 + morse(\o_1^j)} \, + \,  1.
$$
Hence we obtain that
$$
\mu_2 \, = (-1)^{\i_2} + \mu_1 + (-1)^{\i_2} \mu_1 = \,  -1 , 
$$
where we have used \eqref{muu1}.
Therefore the problem $(\mathcal{N}_K)$ has at least one solution in the interval $ \big[ ( 2 S_n (1-\eta ) ) ^{2/n} , ( 2 S_n (1 + \eta ) ) ^{2/n} \big]$. \\
We need to obtain an iteratively  relation between  $ \mu_k$'s. Let $(\o_k^j)$ be the solutions (if there exist) of $(\mathcal{N}_K)$ whose energy levels belong to the interval $[  (k S_n (1- \eta)) ^{2/n},(k S_n (1+ \eta))^{2/n} ]$. For $3 \leq p \leq N$, arguing as above we see that  the level sets $J_{K,\tau}^{  ( p S_n (1+ \eta))^{2/n} }$  and $J_{K,\tau}^{ (p  S_n (1 -\eta))^{2/n} }$ are contractible and blowing up solutions whose energy levels lie in the interval $[  (p S_n (1- \eta)) ^{2/n},(p S_n (1+ \eta))^{2/n} ]$ are:
\begin{itemize}
\item the solutions which converge weakly to a solution $\o_{p-2}$ and blow up at $y_1$ and $y_2$,
\item the solutions which converge weakly to a solution $\o_{p-1}$ and blow up at $y_1$,
\item the solutions which converge weakly to a solution $\o_{p-1}$ and blow up at $y_2$.
\end{itemize}
In addition, there are also the solutions which converge strongly to a solution $\o_p$.
Hence, for $3 \leq p \leq N$, we obtain that
$$
  1= \, \mu_p \, + \, \sum_{\o_{p-2}} (-1)^{2+\i_2 + m(\o_{p-2})} \, + \, \sum_{ \o_{p-1}} (-1)^{1+ m(\o_{p-1})}  \,+  \,  \sum_{ \o_{p-1}} (-1)^{1+ \i_2 + m(\o_{p-1})} \, + \,  1,
$$
which implies 
$$
0 = \mu_p \, - \,  \big( 1+  (-1)^{\i_2}\big) \mu_{p-1} \, + \,  (-1)^{\i_2} \mu_{p-2} \quad \quad \forall \, \, 3 \leq p \leq N . 
$$
Hence two cases may occur:
\begin{description}
\item[(i) ]  $\i_2: = \i(y_2)$ is even. In this case, the last equality implies that $$ \mu_p - \mu_{p-1} = \mu_{p-1} - \mu_{p-2 } = \cdots = \mu_2 - \mu_1 = 0 \quad \forall \, \, 3 \leq p \leq N .
$$ Thus $ \mu_p = \mu_1 = -1$ for each $ 1 \leq p \leq N$. 
\item[(ii) ]  $\i_2: = \i(y_2)$ is odd. In this case, we obtain  $$ \mu_p = \mu_{p-2}  \quad \forall \, \, 3 \leq p \leq N.$$
This implies that, for $ 1 \leq p \leq N$, we have
$$
\mu_p :=
   \begin{cases}
  -1, & \mbox{ if }  p \mbox{ is even} ,   \\
   1 & \mbox{ if } p  \mbox{ is odd}.
\end{cases}
$$
\end{description}
This completes the proof of Theorem \ref{th:tt1}  in the case where $ \#\mathcal{K}^\infty  =2$.
 \end{pfn}

\section{ Counting index formulae and Morse type relations}

Unfortunately  for arbitrarily  $m$ computing $\mu_k$ explicitly   for any number of points in $\mathcal{K}^{\infty}$ is very difficult and to prove that $\mu_k \neq 0$ we have to argue differently using recurrence relations between the $\mu_k$'s. To this aim we make use of the following notation:  For $1 \leq k \leq m:= \# \mathcal{K}^{\infty}$,  we introduce the following  sets:
\be\label{*100}
\mathcal{S}^{\infty, p}_{\geq k}  \qquad; \qquad   \mathcal{S}^{\infty, p}_{\geq k; k}
\ee
\begin{itemize}
\item $ \mathcal{S}^{\infty, p}_{\geq k} $ is the  set of critical points $u_{\infty, p}^\tau$  such that  $J_{K,\tau}( u_{\infty, p}^\tau) \in [( p S_n (1 -\eta))^{2/n}, ( p S_n (1 + \eta)) ^{2/n}]$  and  which blow up at  $r$ points $y_{i_1}, \cdots , y_{i_r}$ with $i_j \geq k$ for all $1 \leq j \leq r$, $r \leq p$ and $u_{\infty, p}^\tau$ converges weakly to $ \o_{p-r}$ (with $\o_0 = 0$),
\item $ \mathcal{S}^{\infty, p}_{\geq k; k} $ is the set of blowing up  critical points belonging to $ \mathcal{S}^{\infty, p}_{\geq k}$  such that $y_k$ is one of their  blow up points.
\end{itemize}
 Furthermore we set
$$  \mu_{\geq k}^{\infty , p } :=
\begin{cases}
  0 , & \mbox{ if }    k > m:= \# \mathcal{K}^\infty ,  \\
  \sum_{ u_{\infty, p}^\tau \in \mathcal{S}_{\geq k} ^ {\infty , p }} (-1)^{ m( u_{\infty, p}^\tau ) }, & \mbox{otherwise} , 
\end{cases}
$$
and
$$ \quad \mu_{\geq k; k} ^{\infty , p } :=
\begin{cases}
  0 , & \mbox{ if }    k > m:= \# \mathcal{K}^\infty , \\
  \sum_{ u_{\infty, p}^\tau \in \mathcal{S}_{\geq k ; k }^ {\infty p }} (-1)^{ m( u_{\infty, p}^\tau ) } , & \mbox{otherwise}.
\end{cases}
$$
We point out that these above quantities are the contribution of the  blowing up critical points in the sets    $ \mathcal{S}^{\infty, p}_{\geq k} $  and $ \mathcal{S}^{\infty, p}_{\geq k; k} $ in the Euler-Poincar\'e characteristic of the pair $$ \big(J_{K,\tau}^{ (p S_n (1 +\eta))^{2/n}}, J_{K,\tau}^{( p S_n (1 -\eta)) ^{2/n} } \big).$$ 
In the next two lemmas, we show   some useful  relations between  $ \mu_p$, $ \mu_{\geq k}^{\infty , p } $ and $  \mu_{\geq k; k} ^{\infty , p } $. Our first relations read as follows:
\begin{lem}[recurrence relations]
Using the  above notation, for each $ 1 \leq k \leq m $, there holds:
\be  \label{1*new} \mu_{\geq k} ^{\infty , 1 } = \sum_{j=k}^m ( -1 )^{\i_j} = ( -1 )^{\i_k} + \mu_{\geq k+1} ^{\infty , 1 } . \ee
Furthermore, for $p \geq 2$, we have 
  \begin{align}
 & \label{*1new} \mu_{\geq k} ^{\infty , p } =  \mu_{\geq k + 1} ^{\infty , p } + \mu_{\geq k; k} ^{\infty , p }  \qquad \mbox{ and } \qquad \mu_{\geq k; k} ^{\infty , p } =  (-1)^{1+ \i(y_k)} \Big[  \mu_{p-1 }+\mu_{\geq k + 1} ^{\infty , p-1 } \Big], \\
 & \label{*1} \mu_{\geq k} ^{\infty , p }  \, = \,\mu_{\geq k + 1} ^{\infty , p }+   (-1)^{1+ \i(y_k)} \Big[  \mu_{p-1 } +\mu_{\geq k + 1} ^{\infty , p-1 }\Big] .
\end{align}
\end{lem}

\begin{pf}
By definition,  a sequence of critical points  $u^{\tau}_{\infty, 1 } \in \mathcal{S}_{\geq k} ^ {\infty , 1 }$ has to blow up at one point $ y_j$ with $ j \geq k$ and has to converge weakly to $0$. Hence the proof of  \eqref{1*new} follows.

Now, for $ p \geq 2$, by definition,  a sequence of critical points  $u^{\tau}_{\infty, p } \in \mathcal{S}_{\geq k} ^ {\infty , p }$ has to blow up and two cases may occur:
\begin{itemize}
\item either $y_k$ does not belong to  the set  of the blow up points of $u^{\tau}_{\infty, p }$. In this case, the function $u^{\tau}_{\infty, p } \in  \mathcal{S}_{\geq k+1} ^ {\infty , p }$ .
 \item or $y_k$  belongs   to the set of the blow up points of $u^{\tau}_{\infty, p }$. In this case, the function $u^{\tau}_{\infty, p } \in  \mathcal{S}_{\geq k; k} ^ {\infty , p }$.
\end{itemize}
Hence  the first claim of \eqref{*1new} follows.\\
Next,  the second above mentioned  situation can be split  into  two subcases:
\begin{itemize}
\item either $y_k$ is the only blow up point of $u^{\tau}_{\infty, p }$ and in this case, $u^{\tau}_{\infty, p }$ has to converge weakly to a critical point $\o_{p-1}$ of $ J_K$.
By Theorem $1$ of \cite{Blow-up sphere} we have that
$$ (-1) ^{\mbox{morse index}(u^{\tau}_{\infty, p }) } = - (- 1)^{\i(y_k)} (-1)^{\mbox{morse index}(\o_{p-1})}.$$
\item or $u^{\tau}_{\infty, p } \rightharpoonup \o_q$ (with $ 0\leq q \leq p-2$)  where  $\o_q$ is either zero (if $q=0$)  or a critical point of $ J_K $ with  $ J_K(\o_q) \in [ ( q S_n (1-\eta)) ^{2/n}, ( q S_n (1+\eta)) ^{2/n}] $ and $u^{\tau}_{\infty, p }$ blows up at $y_k$ and  at some other critical points $y_{j_1}, \cdots , y_{j_{p-q-1}}$ with $j_i \geq k+1$ for each $i$.  Furthermore  according to Theorem $2$ of \cite{Blow-up sphere}, there exists a unique  solution $u^{\tau}_{\infty, p-1 }$ of $ (\mathcal{N}_{K, \tau})$ which blows up at these points $y_{j_1}$, ... , $y_{j_{p-q-1}}$ and converges weakly to $\o_q$ as $\tau \to 0$. Hence $u^{\tau}_{\infty, p-1 }\in \mathcal{S}_{\geq k+1} ^ {\infty , p-1 }$. Hence inferring again to  Theorem   $2$ of \cite{Blow-up sphere}  we have that
$$ (-1) ^{\mbox{morse index}(J_{K,\tau}, u^{\tau}_{\infty, p })} = - (-1)^{\i(y_k)}  (-1) ^{\mbox{morse index}(J_{K,\tau}, u^{\tau}_{\infty, p-1 })}.$$
\end{itemize}
Summing over all  critical points in   $\mathcal{S}_{\geq k} ^ {\infty , p } $ respectively in $\mathcal{S}_{\geq k+1} ^ {\infty , p-1 }$ we recover the second equality in  \eqref{*1new}. Finally we notice that \eqref{*1} follows from the two equalities in \eqref{*1new}.
\end{pf}

The second type of relations  of the above quantities are  some  \emph{Morse equalities} relating the  difference of topology of the pair   $ (J_{K,\tau}^{( p S_n (1 +\eta)) ^{2/n}}, J_{K,\tau}^{ ( p S_n (1 -\eta) )^{2/n} })$  to the topological contribution of all critical points whose energy levels lie in the interval   $[ (p S_n (1 - \eta) )^{2/n}, ( p S_n (1 + \eta) )^{2/n} ]$.  Namely we prove:
\begin{lem}[Morse equalities] With the above notation, we have
\begin{align}
& (i) \quad  \mu _1 +\mu_{\geq 1} ^{\infty , 1 } = 1  \qquad \mbox{ and } \qquad \mu _p+\mu_{\geq 1} ^{\infty , p } = 0 \quad \forall \, \,  2 \leq p \leq  N. \label{*2} \\
&  (ii) \quad  \mu_p +\mu_{\geq 2}^{\infty, p} = 0 \quad \forall \, \,  1 \leq p \leq N.  \label{*3}
  \end{align}
\end{lem}
\begin{pf}
Since $J_{K,\tau}$ satisfies the Palais-Smale condition,  it follows from the classical deformation lemma that
\begin{align*}
& J_{K,\tau}^{ ( p  S_n (1+ \eta) )^{2/n} } \mbox{ retracts by deformation onto } J_{K,\tau}^{ ( p  S_n (1- \eta) )^{2/n} } \bigcup _{u_p^\tau \in \mathcal{S} ^p \cup \mathcal{S}^{\infty, p }_{ \geq 1}} W_u(u_p^\tau)\quad \mbox{ for } p \geq 2, \\
& J_{K,\tau}^{ ( S_n (1+ \eta) )^{2/n} } \mbox{ retracts by deformation onto
} \bigcup _{u_1^\tau \in \mathcal{S}^1 \cup \mathcal{S}^{\infty, 1 }_{ \geq 1}} W_u(u_1^\tau), \end{align*}
where $W_u(u_p^\tau) $ denotes the unstable manifold at the critical point $u_p^\tau$ of $J_{K,\tau}$, $\mathcal{S}^p$ and  $\mathcal{S}^{\infty, p }_{ \geq 1}$ are introduced in \eqref{****} and \eqref{*100} respectively.
Using the Euler-Poincar\'e theorem, we derive that
\begin{align*}
& \chi (  J_{K,\tau}^{ ( p  S_n (1+ \eta) )^{2/n} } ) = \chi (  J_{K,\tau}^{ ( p  S_n (1- \eta) )^{2/n} } ) + \sum_{u_p^\tau \in \mathcal{S}^p \cup \mathcal{S}^{\infty, p }_{ \geq 1}} (-1)^{\mbox{ morse index } (u_p ^\tau)}  \quad \mbox{ for } p \geq 2, \\
& \chi (  J_{K,\tau}^{ ( S_n (1+ \eta) )^{2/n} } ) =  \sum_{u_1^\tau \in \mathcal{S}^1 \cup \mathcal{S}^{\infty, 1 }_{ \geq 1}} (-1)^{\mbox{ morse index } (u_1 ^\tau)}. \end{align*}
Furthermore, Proposition \ref{l:C} implies that $J_{K,\tau}^{ ( p  S_n (1+ \eta) )^{2/n} } $ and $J_{K,\tau}^{ ( p  S_n (1- \eta) )^{2/n} } $ are contractible sets. Therefore, for $p \geq 2$, we get
$$ 1 = 1 + \sum_{u_p^\tau \in \mathcal{S}^p} (-1)^{\mbox{ morse index } (u_p ^\tau)} + \sum_{u_{\infty, p}^\tau \in  \mathcal{S}^{\infty, p }_{ \geq 1}} (-1)^{\mbox{ morse index } (u_{\infty, p} ^\tau)} $$
which gives the second claim in \eqref{*2}. Likewise follows  the first claim of \eqref{*2}. \\
Next  using \eqref{1*new},  \eqref{*1} and \eqref{*2}, we derive that
\be\label{qas2}
(-1)^{\i(y_1)} + \mu_{\geq 2}^{\infty, 1} + \mu_1 = 1 \, \, \,  \mbox{ and } \, \, \,  - (-1)^{\i(y_1)} \Big[ \mu_{\geq 2} ^{\infty , p-1 } + \mu_{p-1 } \Big] + \mu_{\geq 2}^{\infty, p} + \mu_p = 0 \quad \forall \, \, 2 \leq p \leq N.
\ee
Since  $y_1$ is a maximum point of $K$, we have that $ \D K(y_1) <0$ and $ \i (y_1) =0$ is even.  This implies that the first equation in \eqref{qas2} is equivalent to 
$$  \mu_{\geq 2}^{\infty, 1} + \mu_1 = 0.$$ 
Furthermore  denoting  $x_p:=  \mu_{\geq 2}^{\infty, p} + \mu_p$, on one hand the above equation implies that $x_1=0$. On the other hand, we see that the second equation in \eqref{qas2} is equivalent to $ x_p = x_{p-1}$ for each $2 \leq p \leq N$. Thus, we obtain  Claim \eqref{*3} and  the proof of the lemma is thereby complete.
\end{pf}

Next we are going to prove the following crucial Lemmas about the value of $\mu_p$. The first one reads as follows:

\begin{lem}\label{ll:1} 
 Assume that  $ m \geq 3$ and $\i(y_j)$ is even for each $j \geq 2$. Then it holds that 
  $$ \mu _p = - C_{p+m-2}^{m-2}  \quad \mbox{  for each }  p \in \{ 1 , \cdots, N \}.$$
 \end{lem}
 \begin{pf}
 Observe that, since $\i(y_j)$ is even for each $j \geq 2$, from \eqref{1*new} and \eqref{*2}, it holds
 \be \label{mk1} \mu_{\geq k}^{\infty, 1 } = \sum_{j=k} ^m (-1)^{\i(y_j)} = m-k+1 \qquad \mbox{ and } \qquad \mu_1 = 1 - m . \ee 
 Furthermore, \eqref{*1} can be written as 
 $$ \left( \mu_q +\mu_{\geq k}^{\infty, q}  \right) =  \left( \mu_q +\mu_{\geq k+1}^{\infty, q}  \right) -  \left( \mu_{q-1} +\mu_{\geq k+1 }^{\infty, q-1}  \right) \quad \mbox{ for each } q \in \{2, \cdots, N \} $$
 and therefore,  we derive that
\be \label{mk2}  
 \left( \mu_p +\mu_{\geq k+1}^{\infty, p}  \right) = \sum_{q=2}^p \left( \mu_q +\mu_{\geq k}^{\infty, q}  \right) + \left( \mu_1 +\mu_{\geq k+1}^{\infty, 1}  \right)  \qquad \mbox{ for each } p \in \{2, \cdots, N \} .
 \ee
 Now, we claim that 
 \be \label{mk3} 
  \left( \mu_p +\mu_{\geq k+1}^{\infty, p}  \right) = - C_{p+k-2}^{k-2}  \quad \mbox{  for each }  p \in \{ 2 , \cdots, N \} \mbox{ and for each } k \in \{ 2 , \cdots, m \}.
 \ee
 The proof of the claim will be by induction on $k$. \\
 For $ k = 2 $, using \eqref{mk2}, \eqref{*3} and \eqref{mk1}, we get
 $$  \left( \mu_p +\mu_{\geq 3}^{\infty, p}  \right) = \sum_{q=2}^p \left( \mu_q +\mu_{\geq 2}^{\infty, q}  \right) + \left( \mu_1 +\mu_{\geq 3}^{\infty, 1}  \right) = - 1 = - C_p ^0$$
 which proves the claim for $ k = 2 $.\\
 Now, assuming that the claim holds for $ k \in \{ 2 , \cdots, m-1 \}$, then \eqref{mk2} and \eqref{mk1} imply that 
 \begin{align*}  \left( \mu_p +\mu_{\geq k+2}^{\infty, p}  \right)  & = \sum_{q=2}^p \left( \mu_q +\mu_{\geq k+1}^{\infty, q}  \right) + \left( \mu_1 +\mu_{\geq k+2}^{\infty, 1}  \right) \\
 & = - \left( \sum_{q=2}^p C_{q+k-2}^{k-2} \right) - k = - C_{p+(k+1)-2}^{(k+1)-2} .
 \end{align*}
 Hence our claim is proved.\\
 Finally, the proof of Lemma \ref{ll:1} follows from the fact that $ \mu_{\geq m+1}^{\infty, p}  = 0$ by taking $ k=m$ in the previous claim.
\end{pf}

The second tool to determine the value of the local degree $\mu_p$ is provided in the next lemma.

\begin{lem}\label{ll:1bis} 
 Assume that  $ m \geq 3$ and $\i(y_j)$ is odd for each $j \geq 2$. Then it holds that 
  $$ \mu _p =  (-1)^{p + 1} \, C_{p+m-2}^{m-2}  \quad \mbox{  for each }  p \in \{ 2 , \cdots, N \}.$$
 \end{lem}
 \begin{pf}
  Observe that, since $\i(y_j)$ is odd for each $j \geq 2$, from \eqref{1*new} and \eqref{*3}, it holds
 \be \label{mk1bis} \mu_{\geq k}^{\infty, 1 } = \sum_{j=k} ^m (-1)^{\i(y_j)} = - ( m-k+1 ) \qquad \mbox{ and } \qquad \mu_1 = -  \mu_{\geq 2}^{\infty, 1 } =  m - 1 . \ee 
 Furthermore, \eqref{*1} can be written as 
 $$ \left( \mu_q +\mu_{\geq k}^{\infty, q}  \right) =  \left( \mu_q +\mu_{\geq k+1}^{\infty, q}  \right) +   \left( \mu_{q-1} +\mu_{\geq k+1 }^{\infty, q-1}  \right) \quad \mbox{ for each } q \in \{2, \cdots, N \} $$
 and therefore,  we derive that, for each $ k \in \{1, \cdots , m \}$ and  each $ p \in \{2, \cdots, N \}  $, 
\be \label{mk2bis}  
(-1)^p \,  \left( \mu_p +\mu_{\geq k+1}^{\infty, p}  \right) = \sum_{q=2}^p ( -1)^q \, \left( \mu_q +\mu_{\geq k}^{\infty, q}  \right)  -  \left( \mu_1 +\mu_{\geq k+1}^{\infty, 1}  \right)  .
 \ee
 Now, we claim that 
 \be \label{mk3bis} 
 (-1)^p \,  \left( \mu_p +\mu_{\geq k+1}^{\infty, p}  \right) = - C_{p+k-2}^{k-2}  \quad \mbox{  for each }  p \in \{ 2 , \cdots, N \} \mbox{ and for each } k \in \{ 2 , \cdots, m \}.
 \ee
 The proof of the claim will be by induction on $k$. Observe that, for $ k = 2 $, applying \eqref{mk2bis}, we get 
 $$ 
(-1)^p \,  \left( \mu_p +\mu_{\geq 3}^{\infty, p}  \right) = \sum_{q=2}^p ( -1)^q \, \left( \mu_q +\mu_{\geq 2}^{\infty, q}  \right)  -  \left( \mu_1 +\mu_{\geq 3}^{\infty, 1}  \right)  = - 1 = - C_p ^0 , 
 $$
 where we have used \eqref{*3}. Hence the claim \eqref{mk3bis} is true for $ k=2 $. \\
 By induction on $k$, assume that the claim \eqref{mk3bis} holds for $ k-1 \in \{2 , \cdots, m-1 \}$. Applying \eqref{mk2bis} and \eqref{mk1bis},  we obtain 
\begin{align*}
(-1)^p \,  \left( \mu_p +\mu_{\geq k+1}^{\infty, p}  \right) & = \sum_{q=2}^p ( -1)^q \, \left( \mu_q +\mu_{\geq k}^{\infty, q}  \right)  -  \left( \mu_1 +\mu_{\geq k+1}^{\infty, 1}  \right)  \\
 & = \sum_{q=2}^p  \, - C_{ q +k- 3 }^{k-3} - ( k - 1 ) \\
  & =  - \sum_{q=0}^p  \, C_{ q +k- 3 }^{k-3}  = - \,  C_{ p +k- 2 }^{ p } . 
\end{align*}
 Hence, the proof of the claim \eqref{mk3bis} is thereby completed. \\
 Finally, taking $ k = m $ in \eqref{mk3bis} and using the fact that $ \mu_{\geq m+1}^{\infty, p}  = 0 $, we deduce Lemma \ref{ll:1bis}
  \end{pf}
  
  Notice that, in Lemmas \ref{ll:1} and \ref{ll:1bis}, we assumed that, for $ j \geq 2 $, the parameters $ \i (y_j) $'s have the same sign. In the sequel, we need to cover the other situations.

\begin{lem}\label{ll:2}
Assume that   $m \geq 3$, that $\i(y_2)$ is odd and let $ 1 \leq \ell_0\leq (m-1)/2$ be such that $\i(y_{2j})$ is odd and $\i(y_{2j+1})$ is even for $j\in\{1,\cdots, \ell_0\}$. Then there  holds:
$$  \begin{cases} 
\mu_{2p-1}+\mu_{\geq 2j+2}^{\infty, 2p-1} \, = 0  & \forall \, \,  1 \leq p \leq [ (N+1) /2] , \\
 \mu_{2p }+\mu_{\geq 2j+2}^{\infty, 2p}  = - C_{p+j-1} ^p  &  \forall \, \,  1 \leq p \leq [N/2]  , \end{cases}   \qquad \forall \, \,  1 \leq j \leq \ell_0.
$$
\end{lem}

\begin{pf}
First, if $ \i(y_{k})$ is odd and $\i(y_{k+1})$ is even, applying \eqref{*1} two times, it follows that
\begin{align} 
& \left( \mu_q +\mu_{\geq k}^{\infty, q}  \right) = \left( \mu_q +\mu_{\geq k+2}^{\infty, q}  \right) - \left( \mu_{q-2} +\mu_{\geq k+2}^{\infty, q-2}  \right) \qquad \forall \, \, q \in \{ 3, \cdots, N \} , \label{mc1} \\
& \left( \mu_2 +\mu_{\geq k}^{\infty, 2}  \right) = \left( \mu_2 +\mu_{\geq k+2}^{\infty, 2}  \right) + 1 .  \label{mc111}
\end{align}

Second, taking $ k =2 $ in \eqref{mc1} and using \eqref{*1}, we obtain
$$ \left(  \mu_{q} + \mu_{\geq 4}^{\infty, q } \right) =   \left(  \mu_{q-2} + \mu_{\geq 4}^{\infty, q-2} \right) \qquad \forall \, \, q \in \{ 3, \cdots, N \} $$
which implies that (by using \eqref{*2} and \eqref{mc111})
\begin{align}
&  \left(  \mu_{ 2 p -1 } + \mu_{\geq 4 }^{\infty, 2p-1} \right)  =   \left(  \mu_{ 1 } + \mu_{\geq 4 }^{\infty, 1} \right) =  \left(  \mu_{ 1 } + \mu_{\geq 2 }^{\infty, 1} \right)  - (-1)^{\i(y_2)} - (-1)^{\i(y_3)} = 0, \label{mc2} \\
& \left(  \mu_{ 2 p  } + \mu_{\geq 4 }^{\infty, 2p } \right)  =   \left(  \mu_{ 2 } + \mu_{\geq 4 }^{\infty, 2} \right)  = - 1. \label{mc3}
\end{align}
Now, we will prove the lemma by induction on $j$. Notice that, \eqref{mc2} and \eqref{mc3} imply the statements of Lemma \ref{ll:2} in the case $ j = 1 $.
Now, we assume that the result holds for $ j-1 \in \{ 1, \cdots, \ell_0 -1\} $ and  we will prove it for $ j $. Note that $\i(y_{2j})$ is odd and  $\i(y_{2j+1})$ is even. Furthermore, by the induction assumption we have that
 \be\label{*7}  \begin{cases}
\mu_{2p-1}+\mu_{\geq 2j}^{\infty, 2p-1} \, = 0  & \forall \, \,  1 \leq p \leq [ (N+1) /2] , \\
 \mu_{2p }+\mu_{\geq 2j}^{\infty, 2p}  = - C_{p+j-2} ^p  &  \forall \, \,  1 \leq p \leq [N/2]  . \end{cases}  
\ee
Using \eqref{mc1} and the first equation of \eqref{*7}, we deduce that 
$$ \mu_{2p-1} + \mu_{\geq 2j +2 }^{\infty, 2p-1}  = \mu_{2p -3 } + \mu_{\geq 2j + 2 }^{\infty, 2p-3}  \qquad \forall \, \,  2 \leq p \leq [(N+1)/2]  
$$
which implies that, for each $ 2 \leq p \leq [(N+1)/2]  $, 
\be   \label{ab1} \mu_{2p-1} + \mu_{\geq 2j +2 }^{\infty, 2p-1}  = \mu_{ 1 } + \mu_{\geq 2j +2 }^{\infty, 1}  = \mu_{ 1 } + \mu_{\geq 2j  }^{\infty, 1}  - (-1)^{\i(y_{ 2j})} -  (-1)^{\i(y_{ 2j +1 })} = 0   . \ee
In addition, for $ p=1 $,  using \eqref{*3}, we have 
\be \label{ab2} 
\mu_{1} + \mu_{\geq 2j +2 }^{\infty, 1} = \mu_{1} + \mu_{\geq 2 }^{\infty, 1}  - \mu_{\geq 2 }^{\infty, 1} + \mu_{\geq 2j + 2 }^{\infty, 1} 
= \sum_{k = 2 }^{2j + 1 } (-1)^{\i(y_k) } = 0 .\ee
Thus, \eqref{ab1} and \eqref{ab2} complete the proof of the first assertion of the  lemma.\\
 Concerning the second statement, using \eqref{mc1}, we deduce that 
 $$ \left( \mu_{2 q }+\mu_{\geq 2j }^{\infty, 2 q }  \right) = \left( \mu_{2 q }+\mu_{\geq 2j +2 }^{\infty, 2 q } \right)  - \left( \mu_{2 (q-1) }+\mu_{\geq 2j +2 }^{\infty, 2(q-1)}  \right) \qquad \forall \, \,  2 \leq q \leq [ N /2 ]  
$$
which implies that 
  $$ \left( \mu_{2 p }+\mu_{\geq 2j +2 }^{\infty, 2 p }  \right) = \left( \mu_{2  }+\mu_{\geq 2j +2 }^{\infty, 2 } \right)  + \sum_{q=2}^p  \left( \mu_{2 q }+\mu_{\geq 2j }^{\infty, 2 q }  \right) \qquad \forall \, \,  2 \leq p \leq [ N /2 ]   .
$$
 Hence, using \eqref{mc111} and \eqref{*7}, for each $  2 \leq p \leq [ N /2 ]$, we have
  $$ \left( \mu_{2 p }+\mu_{\geq 2j +2 }^{\infty, 2 p }  \right) =  - 1 + \sum_{q=1}^p  \left( \mu_{2 q }+\mu_{\geq 2j }^{\infty, 2 q }  \right) = - 1 -  \sum_{q=1}^p  C_{q+j-2} ^q = -  \sum_{q=0}^p  C_{q+j-2} ^q = - C_{p+j-1} ^p.
$$
Finally, for $ p = 1 $, using \eqref{mc111} and \eqref{*7}, we get 
 $$ \left( \mu_{2  } + \mu_{\geq 2j +2 }^{\infty, 2  }  \right) =  - 1 +  \left( \mu_{2 }+\mu_{\geq 2j }^{\infty, 2  }  \right) = - 1 -  C_{ j - 1 } ^1  = - C_{ j } ^1. $$
This completes the proof of the lemma.
\end{pf}

\section{ Proof of  the theorems in a general case}

This section is devoted to the proof of the theorems in their most generalities. We notice that, when $ m:= \# \mathcal{K}^{\infty}= 2$, we clearly see that $  \text{Index}_K \neq 1$ and therefore we are in the situation of Theorems \ref{th:t2} and \ref{th:t2bis} stated below. These theorems have been proved in Section $3$ in this case. Thus, in the sequel, we assume that  $m:= \# \mathcal{K}^{\infty} \geq 3$.  For arbitrarily  $N \in \N$ and  for $\eta < \frac{1}{2N +1}$ we choose $\e$ as  Proposition \ref{l:C} and assume that $K \in \mathcal{K}_{\e}$ and satisfies the conditions $(H1,H2,H3)$.

We first consider the case  $  \text{Index}_K := \sum_{y\in \mathcal{K}^{\infty}} (-1)^{\i(y)} = 1$, that is, we are going to prove Theorem \ref{th:t1}.

\begin{pfn} {\bf of Theorem \ref{th:t1}}
 Since $  \text{Index}_K =1$, $m $ has to be odd, say $2\ell +1$; $\ell \geq 1$ with $\ell +1$ points satisfying $\i(y)$ even and $\ell$ points satisfying $\i(y)$ odd. Furthermore observe that  by definition  $\mu_{\geq 2\ell+2} ^{\infty , p } = 0$ for each $  1 \leq p \leq N$.\\
Next,  we arrange  the points $y_k$'s of $ \mathcal{K}^{\infty}$ such that $ \i(y_1), \cdots ,  \i(y_{2\ell +1})$ are even and $ \i(y_2),  \cdots ,  \i(y_{2\ell})$ are odd.
We observe that the  assumptions of Lemma \ref{ll:2} are satisfied (by taking $\ell_0 = \ell $). Thus Lemma \ref{ll:2} implies that  
\begin{align*}
& \mu_{\geq 2\ell+2} ^{\infty , 2p-1 } + \mu_{2p-1} =  \mu_{2p-1} = 0  \qquad \forall\, \, \,  1 \leq p \leq [ ( N + 1 ) / 2 ], \\
&   \mu_{\geq 2\ell+2} ^{\infty , 2 p } + \mu_{2p}  =  \mu_{2p } = - C_{p+ \ell - 1 }^p  \quad \forall\, \, \,  1 \leq p \leq [N/2].
\end{align*}
Hence, for each  $ p \in \{ 1, \cdots, [N/2]\}$,  the functional $ J_K $ has at least $ C_{p+ \ell - 1 }^p $ critical points $\o_{2p, k }$'s such that 
$$
J_K(\o_{2p, k })  \in \big[ ( 2 p S_n (1 -\eta) )^{2/n}, ( 2 p S_n (1 + \eta) )^{2/n} \big] \quad \forall\, \, \,  1 \leq k \leq C_{p+ \ell - 1 }^p .
$$
By using \eqref{JKIK1} and \eqref{JKIK2}, it follows that the  Nirenberg problem $(\mathcal{N}_K)$ has at least  $ C_{p+ \ell - 1 }^p $ solutions  whose energies are close to $ \frac{1}{n} 2 p S_n $. \\
Hence Theorem \ref{th:t1} is proved. 
\end{pfn}

From now, we assume that  $  \text{Index}_K  \neq 1$ and without loss of generality, we can assume that $y_1, \cdots , y_m$ are arranged  so that we are in one of the following four cases: \\
 {\bf Case 1.}   all the indices $\i ( y_j)'s$ are even.\\
 {\bf Case 2.}   for $j\geq 2$, all the indices $\i ( y_j)$ are odd.\\
 {\bf Case 3.}   there exists $l\geq 1$ with $2l\leq m-2$ such that 
 $$
 \i ( y_{2j})\,  \mbox{is odd and}\, \i ( y_{2j+1})\,\mbox{ is even for }\,\,  j \in \{1, \cdots ,l\} \, \mbox{ and }\, \i ( y_j) \, \mbox{ is even} \, \,  \forall \,  j \in\{2 \ell +2 ,\cdots, m \}.
 $$
 {\bf Case 4.}    there exists $l\geq 1$ with $2l\leq m-2$ such that 
 $$
 \i ( y_{2j})\,  \mbox{is odd and}\, \i ( y_{2j+1})\,\mbox{ is even for}\,\,   j \in \{1, \cdots ,l\} \, \mbox{ and }\, \i ( y_j) \,\mbox{ is odd for }
j\in\{2l+2,\cdots, m\}, $$
 where $  \i ( y) = n - \mbox{ morse index } (K,y).$
 
 Notice that, in cases $3$ and $4$, we have used the fact that $  \text{Index}_K \neq 1$ which implies that $ m \geq 2 l + 2 $.

To specify the exact number of solutions obtained and their levels, we will present our multiplicity results for each case. Using the above notations, our multiplicity results are as follows:
\begin{thm}\label{th:t2}
Let $n \geq 7$ and $N \in \mathbb{N}$ be an integer. There exists $\varepsilon_N > 0$ such that, if $ 0 < \varepsilon < \varepsilon_N $  and $ K \in \mathcal{K}_\varepsilon $ satisfying   
$\textbf{(H1-H2-H3)}$, $  \text{Index}_K := \sum_{y \in \mathcal{K}^\infty} (-1)^{\i (y)} \neq 1$ and the assumption stated in Case $1$ or Case $2$, then there hold: \\ 
  for any $ p \in\{1,2,...,N\}$, $(\mathcal{N}_K)$ has at least  $ C_{p + m-2} ^{p}$   solutions whose energies are close to $\frac{p}{n}S_n$. \\
 In particular $(\mathcal{N}_K)$ has at least $ C_{ N + m -1 } ^ {m-1 }  - 1 $  solutions.
\end{thm}

\begin{thm}\label{th:t2bis}
Let $n \geq 7$ and $N \in \mathbb{N}$ be an integer. There exists $\varepsilon_N > 0$ such that, if $ 0 < \varepsilon < \varepsilon_N $  and $ K \in \mathcal{K}_\varepsilon $ satisfying  
$\textbf{(H1-H2-H3)}$, $  \text{Index}_K := \sum_{y \in \mathcal{K}^\infty} (-1)^{\i (y)} \neq 1$ and the assumption stated in Case $3$ or Case $4$, then 
\begin{itemize}
\item for each $ p \in \{1, \cdots , [N/2] \}$, $(\mathcal{N}_K)$ hast at least $ C_{p + m - \ell - 2 }^p $ solutions  whose energies are close to $\frac{2p}{n}S_n$, 
\item for each $ p \in \{1, \cdots , [(N+1) / 2] \}$, $(\mathcal{N}_K)$ has at least $ C_{p + m - \ell - 3 }^{p-1} $ solutions  whose energies are close to $\frac{2p-1}{n}S_n$,
\end{itemize}
\end{thm}

Now, we are going to prove Theorems \ref{th:t2} and \ref{th:t2bis}. 

 \begin{pfn}{\bf  \,of Theorem \ref{th:t2}}
 For $ m=2$, the proof in done in Section 3 and for $ m \geq 3$, the proof follows from Lemmas \ref{ll:1} and \ref{ll:1bis}.  
 \end{pfn}

 \begin{pfn}{\bf  \,of Theorem \ref{th:t2bis}}
 We have to consider Cases $3$ and $4$. 
\begin{itemize}
\item In Case $3$, i.e, there exists $\ell \geq 1$ such that $\i(y_{2k+1})$ is even and $\i(y_{2k})$ is odd for each $k \leq \ell$ and $\i(y_j)$ is even for each $ j \geq 2 \ell +2$.

Note that, in this case, the assumptions of Lemma \ref{ll:2} are satisfied for $\ell_0 = \ell$. Hence,  we have
\be\label{*11}  \begin{cases}
\mu_{2p-1}+\mu_{\geq 2 \ell + 2 }^{\infty, 2p-1} \, = 0  & \forall \, \,  1 \leq p \leq [ (N+1) /2] , \\
 \mu_{2p }+\mu_{\geq 2\ell + 2 }^{\infty, 2p}  = - C_{ p+ \ell - 1 } ^p  &  \forall \, \,  1 \leq p \leq [N/2]  . \end{cases}  
\ee
Since $  \text{Index}_K  \neq 1$, we see that $ m \geq 2\ell +2$.  By assumption of the theorem, $ \i(y_k)$ is even for each $ k \geq 2 \ell + 2 $, and therefore, using \eqref{*1}, we get 
\be \label{ab101} 
\left( \mu_{ q } + \mu_{\geq  k }^{\infty, q} \right) = \left( \mu_{ q } + \mu_{\geq  k+1 }^{\infty, q} \right) - \left( \mu_{ q -1} + \mu_{\geq k +1 }^{\infty, q-1} \right)  \, \forall \, \, q \in \{ 2, \cdots, N\}, \forall \, \, k \geq 2 \ell +2 
 \ee
 which implies that  
 \be \label{ab103} 
\left( \mu_{ p } + \mu_{\geq  k+1 }^{\infty, p } \right) =  \sum_{q=2}^p \left( \mu_{ q } + \mu_{\geq  k }^{\infty, q} \right) +  \left( \mu_{ 1 } + \mu_{\geq k+1 }^{\infty, 1} \right)  \, \forall \, \, p \in \{ 2, \cdots, N\}, \forall \, \, k \geq 2 \ell +2 .
 \ee
Notice that, using \eqref{*11} and the fact that $ \i(y_k)$ is even for $ k \geq 2 \ell +2$, we get 
\be \label{ab105} 
 \left( \mu_{ 1 } + \mu_{\geq k+1 }^{\infty, 1} \right) =  \left( \mu_{ 1 } + \mu_{\geq 2 \ell +2 }^{\infty, 1} \right) - \sum _{ j = 2 \ell +2 }^{ k } (-1)^{\i(y_j)} = - ( k - 2 \ell - 1) .
 \ee
Thus, by induction on $ k \in \{ 2 \ell + 2 , \cdots , m \}$ and using \eqref{*11}, we deduce that 
 \be \label{ab100} 
\left( \mu_{ p } + \mu_{\geq  k+1 }^{\infty, p } \right)  < 0  \quad  \, \forall \, \, p \in \{ 1, \cdots, N\}, \forall \, \, k \geq 2 \ell +2 .
 \ee
Now, we need to be more precise. In fact, using \eqref{ab101}, \eqref{ab100} and \eqref{*11}, we deduce that 
 \be \label{ab102} 
\left( \mu_{ 2p +1 } + \mu_{\geq  k+1 }^{\infty, 2 p +1 } \right) \leq   \left( \mu_{ 2 p } + \mu_{\geq  k + 1 }^{\infty, 2 p } \right)   \quad \forall \, \, p \in \{ 1, \cdots, [(N-1) / 2]\}, \forall \, \, k \geq 2 \ell +2 .
 \ee
Thus, taking $ k = m $ in \eqref{ab102} and \eqref{ab105}, it follows that 
  \be \label{ab10212} 
\mu_1 = - (m-2 \ell -1 ) \leq -1 \quad \mbox{ and } \quad  \mu_{ 2p +1 }   \leq   \mu_{ 2 p }  \quad \forall \, \, p \in \{ 1, \cdots, [( N-1 ) / 2 ]\} .
\ee
 
In addition, \eqref{ab103}, \eqref{ab100} and \eqref{*11} imply 
 \be \label{ab106} 
\left( \mu_{ 2p } + \mu_{\geq  k+1 }^{\infty, 2 p  } \right) \leq  \sum_{q=1}^p  \left( \mu_{ 2 q } + \mu_{\geq  k  }^{\infty, 2 q } \right)   - 1 \quad \forall \, \, p \in \{ 1, \cdots, [N/2]\}, \forall \, \, k \geq 2 \ell +2 .
 \ee
Now, we claim that 
 \be \label{ab107} 
\left( \mu_{ 2p } + \mu_{\geq  k+1 }^{\infty, 2 p  } \right) \leq  - C_{p+\ell  + (k-2\ell -2)}^p \quad \forall \, \, p \in \{ 1, \cdots, [N/2]\}, \forall \, \, k \geq 2 \ell +2 .
 \ee
In fact, for $ k = 2 \ell +2 $, using \eqref{ab106} and \eqref{*11}, we deduce that 
 \be \label{ab108} 
\left( \mu_{ 2p } + \mu_{\geq  2\ell + 3 }^{\infty, 2 p  } \right) \leq  \sum_{ q = 1 } ^ p - C_{ q+\ell  - 1 }^q - 1 = - \sum_{ q = 0 } ^ p C_{ q+\ell  - 1 }^q = - C_{ p + \ell  }^ p   \quad \forall \, \, p \in \{ 1, \cdots, [N/2]\} 
 \ee
which implies that Claim \eqref{ab107} is true for $ k = \ell +2 $.\\
By induction on $ k $, assume that Claim \eqref{ab107} is true for $ k $, hence, using \eqref{ab106}, we deduce that, for each  $ p \in \{ 1, \cdots, [N/2]\}  $, 
 \be \label{ab109} 
\left( \mu_{ 2p } + \mu_{\geq  k+2 }^{\infty, 2 p  } \right) \leq  \sum_{ q = 1 } ^ p - C_{ q+\ell  + (k-2\ell -2) }^q - 1 = - \sum_{ q = 0 } ^ p C_{ q+\ell  + (k-2\ell -2) }^q = - C_{ p + \ell + (k + 1 - 2 \ell -2) }^ p  .
 \ee
This completes the proof of Claim \eqref{ab107}.\\
Finally, taking $ k = m $ in \eqref{ab107} and recalling that $ \mu_{\geq m+1} ^{\infty, q } = 0 $ for each $q $, we derive that
 \be \label{ab110} 
 \mu_{ 2p }  = \left( \mu_{ 2p } + \mu_{\geq  m +1 }^{\infty, 2 p  } \right) \leq  - C_{p + m - \ell - 2 }^p \quad \forall \, \, p \in \{ 1, \cdots, [N/2]\}.
 \ee
Hence the proof of Theorem \ref{th:t2bis} follows from \eqref{cardinalp}, \eqref{ab110} and \eqref{ab10212} in the situation where Case 3 holds.

\item In Case $4$, there exists $\ell \geq 1$ such that $\i(y_{2k+1})$ is even and $\i(y_{2k})$ is odd for each $k \leq \ell $ and $\i(y_j)$ is odd for each $ j \geq 2 \ell +2$.
Clearly, in this case, assumptions of Lemma \ref{ll:2} are satisfied for $ \ell_0 = \ell$. Thus, we see that \eqref{*11} holds in this case.
The proof of this case will follow the previous one with some changes because the sign of $ \i(y_k)$ for $ k \geq 2 \ell +2$. Precisely, equations \eqref{ab101} and \eqref{ab103} become
\be \label{ab2101} 
\left( \mu_{ q } + \mu_{\geq  k }^{\infty, q} \right) = \left( \mu_{ q } + \mu_{\geq  k+1 }^{\infty, q} \right) + \left( \mu_{ q -1} + \mu_{\geq k + 1 }^{\infty, q-1} \right)  \, \forall \, \, q \in \{ 2, \cdots, N\}, \forall \, \, k \geq 2 \ell +2 
 \ee
 \be \label{ab2103} 
(-1)^p \left( \mu_{ p } + \mu_{\geq  k+1 }^{\infty, p } \right) =  \sum_{q=2}^p (-1)^q \left( \mu_{ q } + \mu_{\geq  k }^{\infty, q} \right) -\left( \mu_{ 1 } + \mu_{\geq k+1 }^{\infty, 1} \right)  \, \forall \,   2\leq p \leq N, \forall \,  k \geq 2 \ell +2 .
 \ee
Notice that, using \eqref{*11} and the fact that $ \i(y_k)$ is odd for $ k \geq 2 \ell +2$, we get 
\be \label{ab2105} 
 \left( \mu_{ 1 } + \mu_{\geq k+1 }^{\infty, 1} \right) =  \left( \mu_{ 1 } + \mu_{\geq 2 \ell +2 }^{\infty, 1} \right) - \sum _{ j = 2 \ell +2 }^{ k } (-1)^{\i(y_j)} =  ( k - 2 \ell - 1) .
 \ee
Thus, by induction on $ k \in \{ 2 \ell + 2 , \cdots , m \}$ and using \eqref{*11}, we deduce that 
 \be \label{ab2100} 
(-1)^p \left( \mu_{ p } + \mu_{\geq  k+1 }^{\infty, p } \right)  < 0  \quad  \, \forall \, \, p \in \{ 1, \cdots, N\}, \forall \, \, k \geq 2 \ell +2 
 \ee
and therefore, \eqref{ab102} becomes
 \be \label{ab2102} 
- \left( \mu_{ 2p +1 } + \mu_{\geq  k+1 }^{\infty, 2 p +1 } \right) \leq   \left( \mu_{ 2 p } + \mu_{\geq  k + 1 }^{\infty, 2 p } \right)   \quad \forall \, \, p \in \{ 1, \cdots, [( N - 1 ) / 2 ]\}, \forall \, \, k \geq 2 \ell +2 .
 \ee
 Thus, taking $ k = m $ in \eqref{ab2102} and \eqref{ab2105}, it follows that 
  \be \label{ab210212} 
\mu_1 =  (m-2 \ell -1 ) \geq 1 \quad \mbox{ and } \quad  - \mu_{ 2p +1 }   \leq   \mu_{ 2 p }  \quad \forall \, \, p \in \{ 1, \cdots, [( N-1 ) / 2 ]\} .
\ee

In addition, \eqref{ab2103}, \eqref{ab2100} and \eqref{*11} imply 
 \be \label{ab2106} 
\left( \mu_{ 2p } + \mu_{\geq  k+1 }^{\infty, 2 p  } \right) \leq  \sum_{q=1}^p  \left( \mu_{ 2 q } + \mu_{\geq  k  }^{\infty, 2 q } \right)   - 1 \quad \forall \, \, p \in \{ 1, \cdots, [N/2]\}, \forall \, \, k \geq 2 \ell +2 
 \ee
which is exactly the equation \eqref{ab106}.  The sequel of the proof is exactly the same than the previous case (that is Case 3). Thus we will omit it. 
\end{itemize}
The proof of Theorem \ref{th:t2bis} is thereby completed.
\end{pfn}

\end{document}